\documentclass[11pt]{article}%
\usepackage{graphicx}
\usepackage{subcaption}
\usepackage{makeidx}
\usepackage{amssymb}
\usepackage[all]{xy}
\usepackage{float}
\usepackage[latin1]{inputenc}
\usepackage{amsfonts}
\usepackage{amsmath}
\usepackage{amssymb}
\usepackage{url}
\usepackage{hyperref}
\usepackage{graphicx}%
\usepackage{mathabx}
\usepackage{wrapfig}
\usepackage{mathtools}

\usepackage{amsthm}
\usepackage{url}
\usepackage{hyperref}
\providecommand{\U}[1]{\protect\rule{.1in}{.1in}}

\usepackage{xcolor}

\newtheorem{prop}{Proposition}[section]
\newtheorem{cor}[prop]{Corollary}

\newtheorem{lem}[prop]{Lemma}

\newtheorem{theo}[prop]{Theorem}
\newtheorem{rem}[prop]{Remark}

\newcommand{\CC}{\mathbb{C}}

\newcommand{\RR}{\mathbb{R}}

\newcommand{\La}{ {\cal L }}

\newcommand{\Ea}{ {\cal E }}
\newcommand{\Sa}{ {\cal S }}

\newcommand{\Fa}{ {\cal F }}
\newcommand{\Ga}{ {\cal G }}

\newcommand{\Ma}{ {\cal M }}

\newcommand{\Ha}{ {\cal H }}

\newcommand{\point}{\mbox{\LARGE .}}

\newcommand{\cqfd}{\hfill\blbx \\}
\def\blbx{\hbox{\vrule height 5pt width 5pt depth 0pt}\medskip}

\def \RR{\mathbb{R}}

\def \CC{\mathbb{C}}

\begin{document}

  \title{A note on Riccati matrix difference equations }
  \author{P. Del Moral, E. Horton}


\maketitle

\begin{abstract}
Discrete algebraic Riccati equations and their fixed points are well understood and arise in a variety of applications, however, the time-varying equations have not yet been fully explored in the literature. In this article we provide a self-contained study of discrete time Riccati matrix difference equations. In particular, we provide a novel Riccati semigroup duality formula and a new Floquet-type representation for these equations. Due to the aperiodicity of the underlying flow of the solution matrix, conventional Floquet theory does not apply in this setting and thus further analysis is required. We illustrate the impact of these formulae with an explicit description of the solution of time-varying Riccati difference equations and its fundamental-type solution in terms of the fixed point of the equation and an invertible linear matrix map, as well as uniform upper and lower bounds on the Riccati maps. These are the first results of this type for time varying Riccati matrix difference equations.  
\newline

\emph{Keywords} : {\em Riccati matrix difference equations, discrete time algebraic Riccati equation, Sherman-Morrison-Woodbury  inversion identity, Gramian matrix, matrix positive definite maps, Floquet theory, semigroup duality formula, Lyapunov equations.}
\newline

\emph{Mathematics Subject Classification : Primary: 15A24, 15B48, 93B99 ; Secondary: 15A16, 93B25, 93B05.} 

\end{abstract}


\section{Introduction}

Riccati matrix difference  equations are classical in system theory and signal processing, as well as in optimal
control and estimation theory. Their theoretical and numerical analysis is nowadays rather well developed, however it is outside the scope of this note to provide a detailed discussion with an exhaustive list of references regarding these developments and the corresponding application domains. Thus we refer the reader to the review article~\cite{kucera-2} and the seminal books~\cite{hisham,bittanti-2,lancaster} dedicated to the analysis and the applications of continuous and discrete time Riccati equations. 

This article is concerned with the design of a novel semigroup duality relation between discrete time Riccati matrix difference equations (see (\ref{def-Riccati-drift}) and Theorem~\ref{theo-1-intro}).  We also provide a novel Floquet-type normal formulation of a fundamental-type solution associated with these discrete generation Riccati evolution models (see~\eqref{fund-eq} and  and Theorem~\ref{theo-2}).  We illustrate the impact of these formulae with an explicit description of the solution of time-varying Riccati difference equations and its fundamental-type solution in terms of the fixed point of the equation and an invertible linear matrix map (see  Corollary~\ref{cor-1} and Corollary~\ref{cor-2}). We are also able to provide explicit upper and lower bounds on the Riccati map as a consequence of Theorem~\ref{theo-1-intro} and Corollary~\ref{cor-1}. Under some additional and technical invertibility conditions, several direct proofs can be developed using the Sherman-Morrison-Woodbury inversion formula matrix (cf.  section~\ref{sec-prelim}, Remark~\ref{rem-1} and Remark~\ref{rem-2}, as well as the dual-type 
time-reversed difference Riccati models discussed in the article~\cite{ferrante}). 
The analysis of more general models requires one to develop more sophisticated algebraic matrix-inversion techniques (under our assumptions, none of the matrices introduced in (\ref{E-F-def}) are invertible). 

To the best of our knowledge, the semigroup duality relation presented in this article,
including the explicit descriptions of time-varying Riccati difference equations and their fundamental solution in terms of fixed point matrices, are the first results of this type for this class of discrete generation models. Indeed, while the discrete algebraic Riccati equation (DARE) and the associated fixed point are well understood, there has been very little analysis in the time-varying setting. Thus, in this article, we provide a self-contained study of the latter. We also note that the theory developed in this article is crucial to the stability analysis of discrete time Kalman Ensemble filter and thus the results obtained in this article will allow one to develop the theory of such filters, as discussed in~\cite{DEnKF}, to higher dimensions. 

The continuous time version of the Floquet-type formula presented in this article is discussed in the article~\cite{db-21}. We emphasise that the analysis of Riccati difference equations is far more involved than their continuous time counterparts. From a purely mathematical perspective, the algebra of discrete time matrix models is always
more involved. Furthermore, besides some expected lengthy matrix calculations,
several additional inherent difficulties arise in the theoretical analysis of discrete time models. 
For instance, for homogeneous continuous time models, the observability and controllability Gramian matrix functions are invertible for any positive time horizon (under appropriate observability and controllability conditions). However, in the discrete time setting, the invertibility property is only granted for time horizons larger than the dimension of the problem (see for instance the statement and the proof of Proposition~\ref{lu-estimate-Delta-prop}).
In the same vein, the exponential fundamental matrices associated with the  first variational equations in continuous time are invertible for any positive time horizon (see~\cite{db-21}). Conversely, for discrete time models, whenever the drift matrix is not invertible, the Riccati fundamental matrices are not invertible for any time horizon.
 Moreover, in the discrete time setting,
dual Riccati difference equations are not always defined in terms of the inverse of a
Riccati difference equation and may not have a single pair of negative and positive definite fixed points as in the continuous time setting. 

Due to these difficulties, amongst many others, we agree with the comment given in~\cite{bitmead}: {\em ``The (discrete time) Riccati equation is a difficult beast whose behaviour can often be counterintuitive".} \
 Furthermore, we refer to the article~\cite{bitmead} for a popular list of surprisingly false but often admitted conjectures about Riccati difference equations. It is not always simple to find a self-contained, rigorous and easy-to-read study on the regularity and the stability properties of discrete time Riccati matrix equations. On this topic, we refer to the pioneering work of Kalman~\cite{kalman60} and Deyst and Price~\cite{deyst} in the 1960s.
 As noted in \cite{hitz72}, in both of these articles there was a crucial and commonly made error in the proof which invalidated the results, see~\cite{maybeck} and the more recent articles~\cite{rhudy} for a more detailed discussion and references on these issues. This error was repeated  in numerous subsequent works, including in the seminal lecture book of Jazwinski~\cite{jazwinski}. A first correction was noted in a reply to \cite{hitz72}; see also the reply by Kalman \cite{kalman72}.   
 
We refer to the series of articles~\cite{db-21,bucy,wimmer-2,wimmer-3} and the books~\cite{bittanti-2,lancaster}  for a more thorough discussion on Riccati matrix equations.

The rest of the article is set out as follows. In the remainder of this section we will introduce the Riccati matrix differential equations and some of their properties, as well as some notation that will be used throughout. We also introduce the so-called Riccati matrix products (cf.~\eqref{En-def}) that will be fundamental to the stability and regularity analysis of the Riccati difference equations \eqref{def-Riccati-drift}. In section~\ref{statement-sec} we state our main results. As previously mentioned Theorem~\ref{theo-1-intro} provides a duality relation between the Riccati difference equations introduced in~\eqref{def-Riccati-drift} and consequently, we obtain a Lyapunov equation that relates the positive definite fixed points of these equations, given in Corollary~\ref{cor-Lyapunov}. In Theorem~\ref{theo-2} we provide a Floquet-type representation analogous to the continuous time version given in~\cite{db-21}. In section~\ref{sec-prelim} we provide some preliminary results concerning the Riccati maps and the Gramian matrices introduced in~\eqref{Gn-eq}. We also state and prove several useful properties of the Riccati maps that are often used in the literature but for which we have been unable to find a proof. In section~\ref{theo-1-proof}, we prove Theorem~\ref{theo-1-intro} and discuss some simplifications of the proof under slightly stronger conditions on the model parameters. Section~\ref{theo-unif-intro-proof} is concerned with the proof of Corollary~\ref{theo-unif-intro} and finally, the proof of Theorem~\ref{theo-2} is provided in section~\ref{theo-2-proof}.

\subsection{Matrix Differential Riccati Equations}
We denote by $\Ma_{r}=\RR^{r\times r}$ the ring of $(r\times r)$-square matrices with real entries, for some $r\geq 1$, and by $\Ga \La_r\subset \Ma_{r}$ the general linear group of invertible matrices. When there is no chance of confusion, we also slightly abuse notation and denote by $0$ and $I$ the null and identity matrices, respectively, in $\Ma_{r}$ for any dimension $r\geq 1$. We write $A^{\prime}$ to denote the transposition of a matrix $A$. 
 We also let $\Sa_r\subset \Ma_{r}$ denote the subset of symmetric matrices, $\Sa_r^0\subset\Sa_r$ the subset of positive semi-definite matrices, and $\Sa_r^+\subset \Sa_r^0$ the subset of positive definite matrices. We sometimes use the L\" owner partial ordering notation $S_1\geq S_2$ to mean that a symmetric matrix $(S_1-S_2)$ is positive semi-definite (equivalently, $S_2 - S_1$ is negative semi-definite), and $S_1>S_2$ when $(S_1-S_2)$ is positive definite (equivalently, $S_2 - S_1$ is negative definite). Given $S\in \Sa_r^0-\Sa_r^+$ we denote by $S^{1/2}$ a (non-unique) but symmetric square root of $B$ (given by a Cholesky decomposition). When $S\in\Sa_r^+$ we always choose the principal (unique) symmetric square root. 
 
Given some given matrices $(A,R,S)\in (\Ma_{r}\times\Sa^0_r\times \Sa^0_r)$ 
we denote by $\Phi$ and $\widehat{\Phi}$ the matrix monotone maps from $\Sa^0_r$ into itself defined for any $P\in \Sa^0_r$ by
\begin{equation}\label{def-Riccati-drift} 
\Phi(P):=A(I+PS)^{-1}PA^{\prime}+R\quad
\mbox{\rm and}\quad
\widehat{\Phi}(P):=A^{\prime}(I+PR)^{-1}PA+S.
\end{equation}
We denote by $\Phi_n$ and $\widehat{\Phi}_n$ the Riccati evolution semigroups defined by the inductive composition formula 
$\Phi_n=\Phi\circ\Phi_{n-1}$ and  $\widehat{\Phi}_n=\widehat{\Phi}\circ\widehat{\Phi}_{n-1}$, with the convention $\Phi_0=Id=\widehat{\Phi}_0$, the identity map from $\Sa^0_r$ into itself.

In what follows, we assume that 
the  pair $(A,R^{1/2})$ is controllable and $(A,S^{1/2})$ is observable, in the sense that the controllability and observability matrices,
$$
\left[
R^{1/2},
AR^{1/2},
\ldots,A^{r-1}R^{1/2}\right]
\qquad\mbox{\rm and }\qquad
\left[
\begin{array}{c}
S^{1/2}\\
S^{1/2}A\\
\vdots\\
S^{1/2}A^{r-1}
\end{array}\right],
$$ 
have rank $r$. 
Notice that the above rank conditions are trivially met when $R>0$ and $S>0$. When $R$ has the form $R=B\Sigma B^{\prime}$ for some $\Sigma>0$ and some matrix $B$ with appropriate dimensions, the pair $(A,R^{1/2})$ is controllable if and only if $(A,B)$ is controllable.
Also notice that
$(A,S^{1/2})$ is observable if and only if $(A^{\prime},S^{1/2})$ is controllable. Thus, the  pair $(A^{\prime},S^{1/2})$ is controllable and $(A^{\prime},R^{1/2})$ is observable. 
This duality-type relation between the matrices $(A,R,S)$ and $(A^{\prime},S,R)$ is well known in estimation and optimal control theory in the context of linear-Gaussian filtering or linear-quadratic regulation control. The observability condition ensures that
all the coordinates of  a partially observed system state can be recovered from at most $r$-observations. The controllability condition ensures that a controller can control all the directions of a system. 

The Riccati maps presented in (\ref{def-Riccati-drift}) arise in a variety of areas including in optimal control theory and signal processing~\cite{Mikkola, wimmer-4, hisham}. For the convenience of the reader, a brief proof of the symmetry, monotonicity and positive (semi-)definite preserving properties of Riccati maps is provided in section~\ref{sec-prelim}. We also discuss the positive definite preserving property of $\Phi$ in the case where the pair $(A,R^{1/2})$ is controllable. Thus, whenever the  pair $(A,R^{1/2})$ is controllable and $(A,S^{1/2})$ is observable both of Riccati evolution semigroups $\Phi_n$ and $\widehat{\Phi}_n$  are positive definite, in the sense that
\begin{equation}\label{ref-positive}
P>0\Longrightarrow \forall n\geq 1\qquad \Phi_n(P)>0\quad \mbox{\rm and}\quad \widehat{\Phi}_n(P)>0.
\end{equation}
\textcolor{black}{These properties are folklore in the theory of Riccati matrix difference equations however, we have been unable to find any references stating these results.}

We equip the set $\Ma_{r}$ with the spectral norm $\Vert A \Vert=\sqrt{\lambda_{max}(AA^{\prime})}$ where $\lambda_{max}(\cdot)$ denotes the maximal eigenvalue. The minimal eigenvalue is denoted by $\lambda_{min }(\cdot)$.  

We also denote by  $\mbox{\rm Spec}(A)\subset\CC$ the set of eigenvalues of a matrix $A$, and by 
$\rho(A)$ the spectral radius of a matrix $A$ defined by
$$
\rho(A):=\max{\left\{\vert \lambda\vert~:~\lambda\in \mbox{\rm Spec}(A)\right\}}.
$$

It is well-known that the controllability and observability conditions discussed above ensure that the Riccati difference equations (\ref{def-Riccati-drift}) each have a unique positive definite fixed point 
\begin{equation}\label{def-fixed-points}
\Phi(P_{\infty})=P_{\infty}\in \Sa^+_r\quad \mbox{\rm and}\quad
\widehat{\Phi}(\widehat{P}_{\infty})=\widehat{P}_{\infty}\in \Sa^+_r.
\end{equation}
 In addition, the matrices
\begin{equation}\label{def-B}
E:=A(I+P_{\infty}S)^{-1}\quad \mbox{\rm and}\quad \widehat{E}:=A^{\prime}(I+\widehat{P}_{\infty}R)^{-1}
\quad\mbox{\rm satisfy} \quad \rho(E)\vee \rho(\widehat{E})<1.
\end{equation}
The proof of these assertions can be found in any textbook on Riccati equations, see for instance~\cite{kucera,lancaster,kailath} and the more recent book~\cite{hisham}.  In optimal control theory, the matrix $E$ is often called the closed loop-matrix.

\subsection{Algebraic Lyapunov formula and Riccati matrix products}
Consider the matrix map $\Ea$ and $\Fa$ defined for any $P\in \Sa^0_r$ by the formulae
\begin{equation}\label{E-F-def}
\begin{array}{rcl}
\Ea(P):=A(I+PS)^{-1}\in \Ma_r
&\Longrightarrow & E=\Ea(P_{\infty})\in \Ma_r,
\\
&&\\
\Fa(P):=S(I+PS)^{-1}\in \Sa^0_r&\Longrightarrow&  F:=\Fa(P_{\infty})\in \Sa^0_r.~
\end{array}
\end{equation}
In Theorem~\ref{theo-1-intro} we shall see that the fixed points $(P_{\infty},\widehat{P}_{\infty})$ of the Riccati maps $(\Phi,\widehat{\Phi})$ discussed in (\ref{def-fixed-points}) are connected by the discrete time algebraic Lyapunov formula
$$
(P_{\infty}+\widehat{P}_{\infty}^{-1})^{-1}=E^{\prime}(P_{\infty}+\widehat{P}_{\infty}^{-1})^{-1}E+F.
$$
A more general duality-type formula (\ref{key-form}) between the evolution semigroups $(\Phi_n,\widehat{\Phi}_n)$ is also presented in Theorem~\ref{theo-1-intro}.  The above formula provides a way to solve Lyapunov equations 
of the form $E^{\prime}XE+F=X$  with respect to $X$ by computing the fixed point of a dual Riccati equation instead of computing the conventional solution based on series expansions. 

Whenever  $A$ or $S$ is invertible, the matrix $\Ea(P)$ or $\Fa(P)$ is invertible. In this situation, 
a direct proof of the discrete time algebraic Lyapunov formula stated above based on the inversion formula  (\ref{wood}) can be conducted  (see for instance the direct calculations provided in Remark~\ref{rem-1} at the end of section~\ref{theo-1-proof}). In addition, as shown in Remark~\ref{rem-2}, whenever $A$ is invertible
the matrix $P^-_{\infty}:=(-\widehat{P}^{-1}_{\infty})$ is a negative fixed point of $\Phi$.
 In this particular situation, as for continuous time models~\cite{db-21}, the difference
between the positive and negative fixed points $(P_{\infty}-P^-_{\infty})$ solves the the discrete time algebraic Lyapunov formula stated above. 

An alternative proof of this property in terms of time-reversed Riccati difference equations is provided in the article~\cite{ferrante}, which is dedicated to closed forms solutions of the optimal cost and optimal trajectories of a general class of controlled systems based on a judicious parametrisation of all solutions to an extended symplectic system.

We now consider the directed matrix product $\Ea_n(P)$ defined by
\begin{equation}\label{En-def}
\Ea_{n+1}(P):=\Ea_{n}(\Phi(P))\,\Ea(P)\quad\mbox{\rm with}\quad
\Ea_0(P)=I\quad\Longrightarrow \quad \Ea_n(P_{\infty})=E^n.
\end{equation}

Note that this implies that $\Ea_1 = \Ea$. The matrices $\Ea_n(P)$ play a crucial role in the regularity analysis and the stability theory of Riccati difference equations.
For instance, for any $n\geq 0$ and any $P,Q\in \Sa^0_r$ we have the well-known formula
\begin{equation}\label{phi-form}
\Phi_n(P)-\Phi_n(Q)=\Ea_n(P)\,(P-Q)\,\Ea_n(Q)^{\prime},
\end{equation}
whose proof is a simple consequence of (\ref{form-1}). Observe that
 $\Phi_n$ is a smooth matrix functional with a first order Fr\'echet derivative (see~\cite{niclas}) defined for any $P\in\Sa_r^0$ and $H\in\Sa_r$ by
\begin{equation}\label{fund-eq}
\begin{array}{l}
\displaystyle\nabla\Phi_n(P)\cdot H\,=\,\Ea_n(P)\,H\,\Ea_n(P)^{\prime}\\
\\
\Longleftrightarrow\quad \forall n\geq 1\quad
\nabla\Phi_n(P)=\nabla\Phi_{n-1}(\phi(P))\circ \nabla\Phi(P)\quad \mbox{\rm with}\quad
\nabla\Phi_0(P):=Id.
\end{array}
\end{equation}
In the above display, the symbol $``\circ"$ stands for the composition of operators.

From this perspective, the directed matrix product $\Ea_n(P)$ can be seen as  a fundamental solution of
the first variational equation (\ref{fund-eq}) associated with the Riccati difference equation.
Moreover, setting
$$
P_n=\Phi(P_{n-1})\quad \mbox{\rm and}\quad Q=P_{\infty},
$$
the formula \eqref{phi-form} yields the product formula
\begin{equation}\label{phi-form-infty}
\displaystyle P_n-P_{\infty}=\Ea_n(P_0)\,(P_0-P_{\infty})\,(E^n)^{\prime},
\end{equation}
where we may write $\Ea_{n}(P_0)=\Ea(P_{n-1})\ldots  \Ea(P_{1})\, \Ea(P_{0})$.

 The spectral radius $\rho(E)$ is connected to any norm $\Vert\point\Vert$ of the matrix powers $E^n$ arising in (\ref{phi-form-infty}) by
 Gelfand's formula (see for instance~\cite{rota}) given by
\begin{equation}\label{gelfand-form}
\rho(E)=\lim_{k\rightarrow \infty}\Vert E^k\Vert^{1/k}<1,
\end{equation}
where the latter inequality holds due to \eqref{def-B}. Thus, our observability and controllability conditions ensure the exponential decays of  the matrix norms $\Vert E^n\Vert$ towards $0$ for sufficiently large time horizons. 
When  $P_0$ is close to the fixed point $P_{\infty}$, the matrices $\Ea(P_{n})$ are close to 
 $\Ea(P_{\infty})=E$. Hence, the convergence of the directed product $\Ea_n(P_{0})$ towards $0$  depends on the convergence of $P_n$ towards the positive definite fixed point $P_{\infty}$. On the other hand, in view of (\ref{phi-form-infty}), the convergence of $P_n$ towards $P_{\infty}$ also depends on the convergence of the directed product $\Ea_n(P_0)$ to $0$ , as $n\rightarrow\infty$. As a result, to analyse the stability properties of the fundamental matrices discussed in (\ref{fund-eq}), it is crucial to connect more explicitly the directed products  $\Ea_n(P_0)$ to the $n$-power $\Ea_{n}(P_{\infty})=E^n$ of the limiting matrix.

The asymptotic decay rates of $\Ea_n(P_{0})$ towards $0$ can also be discussed in terms of  the generalised spectral radius (a.k.a. the joint spectral radius) of the set of matrices $(\Ea(P_{k}))_{0\leq k<n}$ (see for instance formula (3.3) in~\cite{daubechies}).
 Unfortunately, these extended spectral radius techniques do not apply to our context as they require one to compute the spectral radius or the Lyapunov exponent of the product of any finite subsequence of the unknown sequence of Riccati matrices $(\Ea(P_{n}))_{n\geq 0}$. For a more detailed discussion on the complexity of computing or estimating the extended spectral radius we refer to~\cite{Tsitsiklis}.
 
  Last but not least, since the flow of matrices $n\mapsto \Ea(P_{n})$ is aperiodic as soon as $P_0\neq P_{\infty}$, the conventional Floquet theory, mainly developed for continuous time models~\cite{birttanti,floquet}, cannot be applied nor extended to this type of Riccati matrix difference models (see~\cite{dacuna,sreedhar} for some extensions of Floquet theory to discrete time models). However, in Theorem~\ref{theo-2} we are still able to provide for any time horizon $n\geq r$ a rather surprising Floquet-type normal form of the directed Riccati products:
  $$
 \Ea_{n}(P_0)=E^n~\La_n(P_0)^{-1}\quad\mbox{with some function $\La_n$ \, s.t.}\quad
\sup_{P_0\in \Sa^0_r}\sup_{n\geq r}\Vert \La_n(P_0)^{-1}\Vert<\infty.
  $$
 For a more precise description of the function $\La_n$ we refer the reader to section~\ref{statement-sec} dedicated to the precise statement of our main results.
The above result is an extended version of the Floquet-type formula presented in~\cite{db-21} in the context of continuous time models to discrete Riccati difference equations.

\subsection{Statement of the main results}\label{statement-sec}
Before we present our main results, we first introduce some further notation. We associate with the functions $(\Ea,\Fa)$ introduced in (\ref{E-F-def}) the increasing sequence of Gramian mappings $\Ga_n$ defined sequentially for any $n\geq 1$ and any $P\in \Sa^0_r$ by the recursion
\begin{equation}\label{Gn-def}
\begin{rcases*}
\displaystyle
\Ga_n(P):=\Fa(P)+\Ea(P)^{\prime}\Ga_{n-1}(\Phi(P))~\Ea(P)\in \Sa^0_r,\\
\\
\displaystyle
G_n:=\Ga_n(P_{\infty})\in \Sa^0_r,
\end{rcases*}
\end{equation}
with the initial condition
$$ 
\Ga_0(P)=0\quad\mbox{\rm and} \quad G_0=0 \quad \Longrightarrow \quad \Ga_1(P)=\Fa(P)
\quad\mbox{\rm and} \quad G_1=F,
$$
where the matrix $F$ was introduced in (\ref{E-F-def}).
Note that we may equivalently write
\begin{equation}\label{Gn-eq}
\begin{rcases*}
\displaystyle
\Ga_n(Q)=\sum_{0\leq k<n}\Ea_k(Q)^{\prime}\Fa(\Phi_k(Q))\Ea_k(Q)\in \Sa^0_r,\\
\\
\displaystyle
G_n=E^{\prime}\,G_{n-1}E+F=G_{n-1}+(E^{\prime})^n\,F\,E^n.
\end{rcases*}
\end{equation}
In the above display, the matrices $\Ea_k(Q)$ are the direct products introduced in (\ref{En-def}) and
 $(E,F)$ is the pair of matrices introduced in (\ref{E-F-def}).
 
Finally, we introduce the parallel addition/harmonic-type mean mapping
\begin{equation}\label{def-H}
\begin{rcases*}
\displaystyle
\Ha~:~(P,Q)\in\left(\Sa^+_r\times\Sa^+_r\right)\mapsto
\Ha(P,Q):=(P+Q^{-1})^{-1}\in\Sa^+_r\\
\\
\displaystyle H:=\Ha(P_{\infty},\widehat{P}_{\infty})\in\Sa^+_r.
\end{rcases*}
\end{equation}

Our first main result is a duality-type formula between the Riccati evolution semigroups $(\Phi_n,\widehat{\Phi}_n)$ and an algebraic Lyapunov equation relating the positive definite fixed points $(P_{\infty},\widehat{P}_{\infty})$ introduced in (\ref{def-fixed-points}).
\begin{theo}\label{theo-1-intro}
For any $P,Q\in \Sa^+_r$ and $n\geq 1$ we have the semigroup duality-type formula
\begin{equation}\label{key-form}
\Ha\left(P,\widehat{\Phi}_n(Q)\right)=\Ea_n(P)^{\prime}\,\Ha\left(\Phi_n(P),Q\right)\,
\Ea_n(P)+\Ga_n(P).
\end{equation}
\end{theo}

The proof of the above theorem is provided in section~\ref{theo-1-proof}. Note that from the discussion of the properties of $
\Phi_n$ provided in (\ref{ref-positive}), we have 
$\Phi_n(P) > 0$ and $\widehat{\Phi}_n(P) > 0$ whenever $P > 0$ and thus~\eqref{key-form} is well-defined.

Applying the duality formula (\ref{key-form}) to $(P,Q)=(P_{\infty},\widehat{P}_{\infty})$ and recalling that under our assumptions, the Lyapunov equation (\ref{Lyapunov-eq}) has a unique solution, we check that the matrix  $H$ introduced in (\ref{def-H})
is the unique solution of the 
 Lyapunov equation 
(\ref{Lyapunov-eq}). This yields a direct proof of the following corollary.
\begin{cor}\label{cor-Lyapunov}
The matrix  $H$ introduced in (\ref{def-H})
is the unique solution of the discrete time algebraic Lyapunov equation 
\begin{equation}\label{Lyapunov-eq}
H=E^{\prime}HE+F.
\end{equation}
\end{cor}

Using the Lyapunov fixed point equation (\ref{Lyapunov-eq}), for any $n\geq 0$ we readily check  that
the sequence of Gramian matrices $G_n$ solving the time varying Lyapunov recursion (\ref{Gn-eq}) can alternatively be defined in terms of solution $H$ of the algebraic Lyapunov equation by the formula
\begin{equation}\label{GlH}
G_n=H-(E^n)^{\prime}HE^n.
\end{equation}
Moreover, it also follows that 
\begin{equation}\label{GlH-2}
G_n\leq H = \lim_{n \to \infty}G_n.
\end{equation}

We now formerly define the sequence of linear maps $\La_n$ that were introduced at the end of the previous section,  given for any $P\in  \Sa^0_r$ by
$$
\begin{array}{l}
\displaystyle
\La_n(P):=
I+(P-P_{\infty})G_n\in\Ma_r\\
\\
\displaystyle\Longrightarrow \La_n(P_{\infty})=I
\quad \mbox{and}\quad
\La_n(P)-\La_n(Q)=(P-Q)\,G_n.
\end{array}$$

\begin{theo}(Floquet-type Representation).\label{theo-2}
For any time horizon $n\geq r$, the function $\La_n$ maps $ \Sa^0_r$ into $\Ga \La_r$. In addition,  
for any $P\in \Sa_r^0$ we have the Riccati matrix product formula
\begin{equation}\label{form-Floquet}
\Ea_n(P)=E^n~\La_n(P)^{-1}\quad \mbox{with}\quad \iota:=
\sup_{P\in \Sa^0_r}\sup_{n\geq r}\Vert \La_n(P)^{-1}\Vert<\infty.
\end{equation}
\end{theo}
 The proof of the above theorem is provided in section~\ref{theo-2-proof}. We end this section with some direct consequences of the above results.

Combining (\ref{phi-form}) with (\ref{form-Floquet}), for any time horizon $n\geq r$ and any $P,Q\in \Sa_r^0$  we have
\begin{equation}
\Phi_n(P)-\Phi_n(Q)=E^n~\La_n(P)^{-1}(P-Q)~(\La_n(Q)^{-1})^{\prime}~(E^{n})^{\prime}.
\label{diff-Phi}
\end{equation}
This yields the Lipschitz property
\begin{equation}\label{Lip-ref}
\Vert \Phi_n(P)-\Phi_n(Q)\Vert\leq (\iota\Vert E^n\Vert)^2~\Vert P-Q\Vert.
\end{equation}
Choosing $Q=P_{\infty}$ in~\eqref{diff-Phi}, we obtain the following corollary.
\begin{cor}\label{cor-1}
For any time horizon $n\geq r$ and any $P\in \Sa_r^0$ we have the formula
$$
\Phi_n(P)=P_{\infty}+E^n~\La_n(P)^{-1}(P-P_{\infty})~(E^{n})^{\prime}.
$$
\end{cor}
The above formula  can be seen as an extension of the Bernstein-Prach-Tekinalp formula~\cite{prach-thesis,prach2015infinite} to discrete time Riccati difference equations, see also~\cite{db-21} for the continuous time version of the above result.

Due to (\ref{form-Floquet}), for any $n\geq r$ and $P,Q\in \Sa_r^0$   we also have the product difference formula
\begin{equation}
\Ea_n(P)-\Ea_n(Q)=E^n~
\La_n(P)^{-1}(P-Q)\,G_n\, \La_n(Q)^{-1}.
\label{diff-E}
\end{equation}
This yields the Lipschitz property
\begin{equation}\label{Lip-ref-2}
\Vert \Ea_n(P)-\Ea_n(Q)\Vert\leq \iota^2\Vert E^n\Vert\Vert H\Vert~\Vert P-Q\Vert.
\end{equation}
Again, setting $Q=P_{\infty}$ in~\eqref{diff-E}, we obtain the following corollary.
\begin{cor}\label{cor-2}
For any time horizon $n\geq r$ and any $P\in \Sa_r^0$ we have the formulae
$$
\Ea_n(P)-E^n=E^n~\La_n(P)^{-1}(P-P_{\infty})\,G_n.
$$
\end{cor}

Finally, we provide some surprising uniform estimates of the Riccati semigroup. 

By (\ref{gelfand-form}), for any $\epsilon\in [0,1[$ there exists some parameter  $n_{\epsilon}\geq 1$ such that for any
$n\geq n_{\epsilon}$ we have
\begin{equation}\label{def-n-epsi}
(E^n)^{\prime}\,P_{\infty}^{-1}\,
E^n\leq (1-\epsilon)~ H \widehat{P}_{\infty}^{-1}H.
\end{equation}
With this in mind, the following estimates are a rather straightforward  consequence of Theorem~\ref{theo-1-intro} and Lipschitz property of Riccati maps (\ref{Lip-ref}). 
\begin{cor}\label{theo-unif-intro}
For any $\epsilon\in [0,1[$, any time horizon $m\geq r$ and $n\geq n_{\epsilon}$, as well as for any  $Q\in \Sa^0_r$ we have the uniform estimates 
\begin{equation}\label{theo-unif-intro-form}
G_r\leq \widehat{\Phi}_m(Q)\quad \mbox{and}\quad
\widehat{\Phi}_n(Q)\leq  \epsilon^{-1}\widehat{P}_{\infty}.
\end{equation}
\end{cor}
The proof of the above corollary is provided in section~\ref{theo-unif-intro-proof}. 

As with many of our results, uniform estimates  for continuous time Riccati semigroup are rather well known~\cite{db-17} however, to the best of our knowledge the estimates presented in (\ref{theo-unif-intro-form}) are completely new for discrete time Riccati equations. In continuous time, these uniform estimates are obtained by sophisticated Riccati differential equation comparisons involving Gramian inversion techniques. 
In this case, the upper bound follows fairly easily from the duality formula given in Theorem~\ref{theo-1-intro}.

\section{Some preliminary results}\label{sec-prelim}
We recall the celebrated Sherman-Morrison-Woodbury matrix sum inversion identity
\begin{equation}\label{wood}
(M+UNV)^{-1}=M^{-1}-M^{-1}U(N^{-1}+VM^{-1}U)^{-1}VM^{-1},
\end{equation}
which is valid for any invertible matrices $(M,N)$ and any  conformable matrices $(U,V)$, see for instance the seminal articles~\cite{bartlett,sherman,woodbury}, an earlier work by Guttmann~\cite{guttman} and the review articles~\cite{henderson,hager}. Several extensions of the above formula to Hilbert state spaces in terms of Moore-Penrose or generalised Drazin inverses can also be found in the article~\cite{deng}. 
Also recall that the eigenvalues of the product $PQ$ of positive semi definite matrices 
$P,Q\in \Sa^0_r$ are nonnegative, so those of $(I+PQ)$ are positive. This elementary property ensures that $(I+PQ)$ is invertible. Thus, the Riccati maps (\ref{def-Riccati-drift})
are well defined without appealing to Moore-Penrose or other types of generalised inverses.

\begin{lem}\label{lem-alpha}
For any $P\in\Sa_r^0$ we have 
$$
\alpha_-(P)S\leq \Fa(P)=S^{1/2}(I+S^{1/2}PS^{1/2})^{-1}S^{1/2}\leq \alpha_+(P)S,
$$
with the positive parameters
$$
\alpha_-(P):=(1+\lambda_{\tiny max}(P)\lambda_{\tiny max}(S))^{-1}~\quad \mbox{and}\quad
\alpha_+(P):=(1+\lambda_{\tiny min}(P)\lambda_{\tiny min}(S))^{-1}.
$$
\end{lem}
\proof
Applying (\ref{wood}) with $M=N=I$, $U=S^{1/2}$ and $V=PS^{1/2}$ we obtain
$$
(I+S^{1/2}PS^{1/2})^{-1}=I-S^{1/2}(I+PS)^{-1}PS^{1/2},
$$ 
and therefore
$$
S^{1/2}(I+S^{1/2}PS^{1/2})^{-1}S^{1/2}=S\left(I-(I+PS)^{-1}PS\right).
$$
Now applying (\ref{wood}) with $M=N=U=I$ and $V=PS$ we also check that
$$
 (I + PS)^{-1} =I-(I+PS)^{-1}PS.
$$ 
We conclude that
\begin{equation}\label{ref-S-P-2}
\Fa(P) = S(I + PS)^{-1} = S\left(I-(I+PS)^{-1}PS\right) = S^{1/2}(I+S^{1/2}PS^{1/2})^{-1}S^{1/2}.
\end{equation}
Now note that
$$
(1+\lambda_{\tiny min}(P)\lambda_{\tiny min}(S))~I\leq 
I+\lambda_{\tiny min}(P)~S\leq I+S^{1/2}PS^{1/2}.
$$
In the same vein, we have
$$
 I+S^{1/2}PS^{1/2}\leq I+\lambda_{\tiny max}(P)~S\leq
(1+\lambda_{\tiny max}(P)\lambda_{\tiny max}(S))~I.
$$
This ends the proof of the lemma.
\cqfd

\begin{lem}
For any $P,Q\in \Sa^0_r$ we have the formulae
\begin{equation}\label{form-1}
\Ea(Q)=\Ea(P)~\left(I+(P-Q)\Fa(Q)\right)
\quad
\mbox{\rm and}
\quad
\Phi(P)-\Phi(Q)=\Ea(P)(P-Q)\Ea(Q)^{\prime}.
\end{equation}
\end{lem}
\proof
The first assertion is a direct consequence of the formulae
$$
\Ea(Q)=\Ea(P) (I+PS)(I+QS)^{-1}
\quad\mbox{\rm and}\quad
(I+PS)(I+QS)^{-1}=\left(I+(P-Q)\Fa(Q)\right).
$$
The second assertion can be verified using the formulae
\begin{eqnarray*}
(I+PS)^{-1}P&=&P(I+SP)^{-1}\\
(I+PS)^{-1}P-Q\,(I+SQ)^{-1}&=&(I+PS)^{-1}(P-Q)~(I+SQ)^{-1}.
\end{eqnarray*} \cqfd

\begin{prop}\label{lu-estimate-Delta-prop}
For any time horizon $n\geq r$ we have
\begin{equation}\label{lu-estimate-Delta}
0< G_r\leq G_n \leq H
 \end{equation}
\end{prop}
\proof
Thanks to (\ref{GlH}) and  Lemma~\ref{lem-alpha}, for any $n\geq r$ we have
$$
\alpha_-(P_{\infty})~\Omega_r\leq
\alpha_-(P_{\infty})~\Omega_n~\leq G_n\leq H,
$$
with the Gramian matrix
$$
\Omega_n:=\sum_{0\leq k<n}(E^k)^{\prime}SE^k\leq \Omega:=\sum_{n\geq 0 }(E^n)^{\prime}SE^n.
$$
Applying (\ref{wood}) to $M=N=I$ and $(U,V)=(P_{\infty},S)$ we have
$$
(I+P_{\infty}S)^{-1}=I-(P_{\infty}^{-1}+S)^{-1}S,
$$
which yields the formula
$$
E=A-A(P_{\infty}^{-1}+S)^{-1}S.
$$
Thus, for any $z\in \CC^r$ and $\lambda\in\CC$ we have
$$
 Ez=\lambda z\quad \mbox{\rm and}\quad S^{1/2}z=0\Longleftrightarrow
  Az=\lambda z\quad \mbox{\rm and}\quad S^{1/2}z=0\Longrightarrow z=0.
$$
Recalling that $(A,S^{1/2})$ is observable, the above equivalence is a direct consequence of 
the Popov-Belevitch-Hautus observability test,~\cite{camlibel}. This ensures that
 the pair $(E,S^{1/2})$ is also observable, that is 
$$
\mbox{\rm Rank}\left[
S^{1/2},
S^{1/2}E,
\ldots,S^{1/2}E^{r-1}\right]=r.
$$
This rank condition ensures that
$
\Omega_r>0$.
Indeed, we have
$$
\begin{array}{l}
x^{\prime}\Omega_rx=0\\
\\
\Longrightarrow  (S^{1/2}x)^{\prime} (S^{1/2}x)=0,~(S^{1/2}Ex)^{\prime}(S^{1/2}Ex)=0,\ldots,
(S^{1/2}E^{r-1}x)^{\prime}(S^{1/2}E^{r-1}x)=0\\
\\
\Longrightarrow S^{1/2}x=0,~S^{1/2}Ex=0,\ldots,~
S^{1/2}E^{r-1}x=0\Longrightarrow x=0
\end{array}
$$
This implies that
$$
\forall n\geq  r\qquad \Omega_n>\Omega_r>0\quad\mbox{\rm and}\quad G_r
\geq \Omega_r^-:= \alpha_-(P_{\infty})~\Omega_r>0.
$$
This ends the proof of the lower bound estimate stated in the left hand side of (\ref{lu-estimate-Delta}) and thus the proposition. \cqfd

For completeness and for the convenience of the reader, to end this section we prove some rather well-known properties of the map $\Phi$, starting with the following lemma.

\begin{lem}\label{lem-phi-props}
The Riccati map $\Phi$ introduced in~\eqref{def-Riccati-drift} satisfies the following properties.
\begin{enumerate}
\item[(i)] For all $P \in \Sa_r^0$, $\Phi(P)' = \Phi(P)$. 
\item[(ii)] For $P, Q \in \Sa_r^0$ and $n \ge 1$, $P \ge Q \Longrightarrow \Phi_n(P) \ge \Phi_n(Q)$ and $\Phi_{n+1}(0) \ge \Phi_n(0)$.
\item[(iii)] For all $P \in \Sa_r^0$, $\Phi(P) \ge 0$. If, in addition, $R$ is invertible, then $\Phi(P) > 0$.
\end{enumerate}
\end{lem}

\proof
\begin{enumerate}
\item[(i)] Applying (\ref{wood}) to  $M=I$, $N=I$, $U=PS^{1/2}$ and $V=S^{1/2}$ we check that
$$
(I+PS)^{-1}P=P-PS^{1/2}(I+S^{1/2}PS^{1/2})^{-1}S^{1/2}P=P(I+SP)^{-1},
$$
from which it follows that $\Phi(P)' = \Phi(P)$.

\item[(ii)] First note that for any $P \in \Sa_r^0$ and any $A \in \Ma_r$ the product $A'PA$ is positive semi-definite. 

By Lemma~\ref{lem-alpha} and equation (\ref{form-1}) we have
\begin{eqnarray*}
  \Phi(P) - \Phi(Q) 
 &=& \Ea(P)(P-Q)\Ea(Q)^{\prime} \\
 &=&\Ea(P)(P-Q)\left(\left(I+\Fa(Q)^{\prime}(P-Q)\right)\right)\Ea(P)^{\prime}\\
 &=&\Ea(P)(P-Q)\Ea(P)^{\prime}+\Ea(P)(P-Q)\Fa(Q)^{\prime}(P-Q)\Ea(P)^{\prime}.
\end{eqnarray*}
{Since $\Fa(Q)$ is positive semi-definite the preceding comments imply that
\begin{equation}\label{red-monotone}
P\geq Q\Longrightarrow  \Phi(P) \geq \Phi(Q) \Longrightarrow\forall n\geq 1\quad
\Phi_n(P) \geq \Phi_n(Q).
\end{equation}

Using (\ref{red-monotone}) we readily check by induction that $\Phi_n(0)$ is a non-decreasing sequence; that is for any $n\geq 1$, we have 
$$
0\leq R\leq   \Phi_n(0)\leq \Phi_{n+1}(0).
$$
} 

\item[(iii)] From the monotonicity of the map $\Phi$, we have
$$
P\geq 0\Longrightarrow \Phi(P)\geq \Phi(0)=R\geq 0.
$$
In addition, if $R>0$, this inequality is strict.
\end{enumerate} \cqfd

We note that even though $\Phi$ is always monotone, the same cannot be said of the sequence $P_n$ unless $P_0=0$.

Next, consider the positive definite preserving properties of $\Phi$ under the condition that  the pair $(A,R^{1/2})$ is controllable. We start with a technical lemma that is interesting in its own right.
We emphasise that this result should be known but we have not been able to find it in the literature.
\begin{lem}\label{lem-A-R}
Whenever the pair $(A,R^{1/2})$ is controllable we have the following property
 $$
 AA^{\prime}+R>0.
 $$
\end{lem}
\proof
Consider a sequence $W_n$ of independent centered Gaussian random variables on $\RR^r$ with unit variance, and set
$$
X_n:=AX_{n-1}+R^{1/2}W_n=A^nX_0+\sum_{0\leq k<n}A^kR^{1/2}W_{n-k}.
$$
In the above display, $X_0$ stands for a centered Gaussian random variable on $\RR^r$ with covariance $P_0\in \Sa_r^0$. We also let $P_n$ denote the covariance of the random variables $X_n$. In this notation,
 the controllability condition ensures that $P_n$ is invertible for any $n\geq r$. Equivalently, we have that
$$
\forall n\geq r\quad\forall P_0\in \Sa^0_r\qquad \lambda_{\tiny min}(P_r)>0.
$$
The covariance $P_n$ of the random variables $X_n$ also satisfies for any $n\geq 1$ the recursion
$$
P_n=AP_{n-1}A^{\prime}+R.
$$
This implies that
$$
0<\lambda_{\tiny min}(P_r)~ I\leq 
P_{r}\leq \lambda_{\tiny max}(P_{r-1}) ~AA^{\prime}+R
$$
from which we conclude that $AA^{\prime}+R>0$.
This ends the proof of the lemma.\cqfd

We are now in a position to prove (\ref{ref-positive}). Whenever $P>0$ we have
$$
\Phi(P)=A\left(P^{-1}+S\right)^{-1}A^{\prime}+R\geq 
\lambda_{\tiny min}(P) ~A\left(I+\lambda_{\tiny min}(P) S\right)^{-1}A^{\prime}+R.
$$
If the pair $(A,R^{1/2})$ is controllable, Lemma~\ref{lem-A-R} ensures that
$$
P>0\Longrightarrow
\Phi(P)\geq \frac{\lambda_{\tiny min}(P)}{1+\lambda_{\tiny min}(P)\lambda_{\tiny max}(S)} ~AA^{\prime}+R>0\Longrightarrow\forall n\geq 1\quad \Phi_n(P)>0.
$$
Analogous arguments clearly show that this also holds for $\widehat{\Phi}_n$. This ends the proof of (\ref{ref-positive}).
\cqfd

\section{Duality-type formulae}\label{theo-1-proof}

Here we provide the proof of Theorem~\ref{theo-1-intro}, followed by some comments on certain, albeit stronger, conditions that greatly simplify the proof of the Lyapunov equation (\ref{Lyapunov-eq}). 

\proof We use an induction argument with respect to the parameter $n$.
Firstly, we check that
the duality formula (\ref{key-form}) is satisfied for $n=1$. Observe that for any $P,Q\in \Sa^0_r$ we have
$$
I+\,P\,\widehat{\Phi}(Q)\,=\left(I+PS\right)+PA^{\prime}Q(I+RQ)^{-1}A.
$$
Applying (\ref{wood}) with $M=\left(I+PS\right)$, $N=(I+RQ)^{-1}$, $U=PA^{\prime}Q$ and $V=A$ we obtain
\begin{equation}\label{rk-1}
\begin{array}{l}
(I+\,P\,\widehat{\Phi}(Q))^{-1}\\
\\
=\left(I+PS\right)^{-1}-\left(I+PS\right)^{-1}PA^{\prime}Q\left((I+RQ)+A\left(I+PS\right)^{-1}PA^{\prime}Q\right)^{-1}A\left(I+PS\right)^{-1}\\
\\
=\left(I+PS\right)^{-1}-P~\Ea(P)^{\prime}Q\left(I+\Phi(P)Q\right)^{-1}\Ea(P).
\end{array}\end{equation}
The last assertion comes from the fact that
$$
A\left(I+PS\right)^{-1}PA^{\prime}=\Phi(P)-R\quad \mbox{\rm and}\quad
\left(I+PS\right)^{-1}P=P\left(I+SP\right)^{-1}.
$$
On the other hand, we have
$$
(I+\,P\,\widehat{\Phi}(Q))^{-1}=I-P\,\widehat{\Phi}(Q)(I+\,P\,\widehat{\Phi}(Q))^{-1}\quad \text{and}\quad
\left(I+PS\right)^{-1}=I-PS\left(I+PS\right)^{-1}.
$$
Using (\ref{rk-1}) we obtain the formula
$$
P~\widehat{\Phi}(Q)(I+\,P\widehat{\Phi}(Q))^{-1}
=P~\Fa(P)+P~\Ea(P)^{\prime}Q\left(I+\Phi(P)Q\right)^{-1}\Ea(P).
$$
When $P\in\Sa^+_r$ this implies that for any $Q\in \Sa^0_r$ we have
$$
\widehat{\Phi}(Q)(I+\,P\widehat{\Phi}(Q))^{-1}
=\Fa(P)+\Ea(P)^{\prime}Q\left(I+\Phi(P)Q\right)^{-1}\Ea(P).
$$
Using (\ref{ref-positive}) for any $Q>0$ we check that
$$
Q\left(I+\Phi(P)Q\right)^{-1}=\Ha\left(\Phi(P),Q\right)\quad
\mbox{\rm and}\quad
\widehat{\Phi}(Q)(I+\,P\widehat{\Phi}(Q))^{-1}=\Ha\left(P,\widehat{\Phi}(Q)\right).
$$
Recalling that $\Ga_1=\Fa$, this ends the proof of the duality formula (\ref{key-form}) for $n=1$.

Now suppose that the duality formula (\ref{key-form}) holds for some $n \ge 1$. Then, we have
\begin{equation*}
\Ha\left(P,\widehat{\Phi}_{n+1}(Q)\right)=\Ha\left(P,\widehat{\Phi}_{n}(\widehat{\Phi}(Q))\right)=\Ea_n(P)^{\prime}\,\Ha\left(\Phi_n(P),\widehat{\Phi}(Q)\right)\,
\Ea_n(P)+\Ga_n(P).
\end{equation*}
On the other hand, applying (\ref{key-form}) to $n=1$ and recalling that $\Ga_1=\Fa$,
we have
$$
\Ha\left(\Phi_n(P),\widehat{\Phi}(Q)\right)=
\Ea(\Phi_n(P))^{\prime}\,\Ha\left(\Phi_{n+1}(P),Q\right)\,
\Ea(\Phi_n(P))+\Fa(\Phi_n(P)).
$$
This implies that
\begin{eqnarray*}
\Ha\left(P,\widehat{\Phi}_{n+1}(Q)\right)&=&\Ea_{n+1}(P)^{\prime}\,\Ha\left(\Phi_{n+1}(P),Q\right)\,
\Ea_{n+1}(P)\\
&&\hskip3cm+\left(\Ga_n(P)+\Ea_n(P)^{\prime}\,\Fa(\Phi_n(P))
\,
\Ea_n(P)\right)\\
&=&\Ea_{n+1}(P)^{\prime}\,\Ha\left(\Phi_{n+1}(P),Q\right)\,
\Ea_{n+1}(P)+\Ga_{n+1}(P).
\end{eqnarray*}
This shows that formula (\ref{key-form}) is valid at rank $(n+1)$, thus for any $n\geq 1$. 
\cqfd

\begin{rem} \label{rem-1}
Whenever  $A$ or $S$ is invertible, a direct proof of (\ref{Lyapunov-eq}) based on the inversion formula  (\ref{wood}) can be conducted.
For instance, when $A$ is invertible, using the fixed point equations, we check that
\begin{equation}\label{ref-A-inv}
\begin{array}{l}
(P_{\infty}^{-1}+S)^{-1}+(\widehat{P}_{\infty}-S)^{-1}\\
\\
\quad=A^{-1}(P_{\infty}-R)(A^{\prime})^{-1}+
A^{-1}(\widehat{P}_{\infty}^{-1}+R) (A^{\prime})^{-1}=A^{-1}H^{-1}(A^{\prime})^{-1}.
\end{array}
\end{equation}
On the other hand, applying (\ref{wood}) with $M=(P_{\infty}^{-1}+S)$, $N=(\widehat{P}_{\infty}-S)^{-1}$
and $U=I=V$ we check that
$$
\left((P_{\infty}^{-1}+S)^{-1}+(\widehat{P}_{\infty}-S)^{-1}\right)^{-1}=(P_{\infty}^{-1}+S)-(P_{\infty}^{-1}+S)~H(P_{\infty}^{-1}+S).
$$
In the same vein, applying (\ref{wood}) with $M=P_{\infty}^{-1}$, $N=\widehat{P}_{\infty}$
and $U=I=V$ we check that
$$
\begin{array}{l}
(P_{\infty}^{-1}+\widehat{P}_{\infty})^{-1}
=P_{\infty}\left(P_{\infty}^{-1}-H\right)P_{\infty}.
\end{array}
$$
Inverting (\ref{ref-A-inv}), this yields the formula
\begin{eqnarray*}
A^{\prime}HA&=&(P_{\infty}^{-1}+S)-(I+SP_{\infty})~
\left(P_{\infty}^{-1}-H\right)(I+P_{\infty}S).
\end{eqnarray*}
This implies that
$$
A^{\prime}HA+(I+SP_{\infty})S=(I+SP_{\infty})~
H(I+P_{\infty}S).
$$
from which we readily check that $H
$ solves the Lyapunov equation (\ref{Lyapunov-eq}).

In the same vein, when $S>0$ is invertible, we have
$$
\widehat{P}_{\infty}=S+A^{\prime}(\widehat{P}_{\infty}^{-1}+R)^{-1}A=S+A^{\prime}(H^{-1}-(P_{\infty}-R))^{-1}A.
$$
Applying (\ref{wood}) with $M=S$, $N=(H^{-1}-(P_{\infty}-R))^{-1}$,
 $U=A^{\prime}$, and $A=V$ we check that
\begin{eqnarray*}
\widehat{P}_{\infty}^{-1}
&=&S^{-1}-S^{-1}A^{\prime}\left(H^{-1}+A\left(S^{-1}-(P_{\infty}^{-1}+S)^{-1}\right)A^{\prime}\right)^{-1}AS^{-1}.
\end{eqnarray*}
On the other hand, applying (\ref{wood}) with $M=P_{\infty}^{-1}$, $N=S$, and
 $U=I=V$  we have
$$
(P_{\infty}^{-1}+S)^{-1}=S^{-1}-S^{-1}(P_{\infty}+S^{-1})^{-1}S^{-1}.
$$
This yields the formula
\begin{eqnarray*}
H^{-1}
&=&(P_{\infty}+S^{-1})-(S^{-1}A^{\prime})\left(H^{-1}+(AS^{-1})(P_{\infty}+S^{-1})^{-1}(S^{-1}A^{\prime})\right)^{-1}(AS^{-1}).
\end{eqnarray*}
We check that $H$ solves the Lyapunov equation 
(\ref{Lyapunov-eq}) by
 applying (\ref{wood}) to the collection of matrices
$$
M=(P_{\infty}+S^{-1})^{-1}\qquad U^{\prime}=V=(A(I+P_{\infty}S)^{-1})\quad \mbox{\rm and}\quad N=H^{-1}.
$$

 \end{rem}
 
 \begin{rem}\label{rem-2}

The map $\Phi$ can be extended to invertible matrices $P\in \Ga \La_r\cap\Sa_r$ s.t. $P+S\in \Ga \La_r\cap\Sa_r$ by setting
$$
\Phi(P):=A(P^{-1}+S)^{-1}A^{\prime}+R.
$$
If we further assume that  $A$ is invertible, we have
\begin{eqnarray*}
\widehat{P}_{\infty}^{-1}+R>0&\Longrightarrow&
\widehat{P}_{\infty}-S:=A^{\prime}(\widehat{P}_{\infty}^{-1}+R)^{-1}A>0 \\
&\Longrightarrow&
\widehat{P}_{\infty}A(\widehat{P}_{\infty}-S)^{-1}A^{\prime}
=\widehat{P}_{\infty}A\left(A^{\prime}(\widehat{P}_{\infty}^{-1}+R)^{-1}A\right)^{-1}A^{\prime}=I+\widehat{P}_{\infty}R\\
&\Longrightarrow&A(\widehat{P}_{\infty}-S)^{-1}A^{\prime}
=\widehat{P}_{\infty}^{-1}+R>0,
\end{eqnarray*}
from which we readily check that
$$
-S^{-1}<{P}_{\infty}^-:=-\widehat{P}_{\infty}^{-1}<0\Longrightarrow 
{P}_{\infty}^-=A\left(({P}_{\infty}^-)^{-1}+S\right)^{-1}A^{\prime}+R
=\Phi({P}_{\infty}^-)<0.
$$
This shows that ${P}_{\infty}^-$ is a negative definite solution of the fixed point equation
$\Phi(P)=P$. Whenever $A$ is invertible, the matrix $(-{P}_{\infty}^-)$  can also be interpreted as the positive definite solution of a time-reversed Riccati difference equation which can be interpreted as a dual Riccati equation~\cite{ferrante}. 
\end{rem}
\section{Uniform estimates}\label{theo-unif-intro-proof}

We now prove the uniform estimates stated in Corollary~\ref{theo-unif-intro}.

\textcolor{black}{
The dual Riccati semigroup $\widehat{\Phi}_n$ is defined as $\Phi_n$ by replacing
the matrices $(A,R,S)$ by $(A^{\prime},S,R)$. In this section, we use the notation $\widehat{\point}$ to the denote the dual mathematical objects; for instance
 $(\widehat{\iota},\widehat{E},\widehat{G}_n)$ stands for the dual parameters
 defined as $(\iota,E,G_n)$ by replacing
$(A,R,S)$ by $(A^{\prime},S,R)$.}

Using Theorem~\ref{theo-1-intro} and Proposition~\ref{lu-estimate-Delta-prop}, for any $Q \in \Sa_r^+$ and $n\geq r$ we have
\begin{equation*}
  \widehat{\Phi}_n(Q) > \Ha\left(P_\infty,   \widehat{\Phi}_n(Q)\right) = (E^n)^{\prime}\,\Ha\left(P_{\infty},Q\right) E^n+G_n > G_r.
\end{equation*}
For any $\epsilon>0$, $n\geq r$ and $Q\in \Sa^0_r$ this implies that
$$
  \widehat{\Phi}_n(Q+\epsilon I)>G_r.
$$
On the other hand, using the Lipschitz property
(\ref{Lip-ref}) we have
$$
\Vert    \widehat{\Phi}_n(Q+\epsilon I)-   \widehat{\Phi}_n(Q)\Vert\leq (\widehat{\iota}~\Vert \widehat{E}^n\Vert)^2~\epsilon,
$$
which yields the uniform estimate
$$
G_r<  \widehat{\Phi}_n(Q)+(\widehat{\iota}~\Vert \widehat{E}^n\Vert)^2~\epsilon I\longrightarrow_{\epsilon\rightarrow 0} \widehat{\Phi}_n(Q).
$$
This ends the proof of the lower bound estimate stated in (\ref{theo-unif-intro-form}).

For the upper bound, applying (\ref{key-form}) to $P=P_{\infty}$ for any $Q>0$ we check that
\begin{equation}
\Ha\left(P_{\infty},\widehat{\Phi}_n(Q)\right)=(E^n)^{\prime}\,\Ha\left(P_{\infty},Q\right)\,
E^n+G_n.
\label{eq-dual}
\end{equation}
On the other hand, by (\ref{GlH}) for any $n\geq 1$ we have
$$
G_n\leq H=(P_{\infty}+\widehat{P}_{\infty}^{-1})^{-1}\quad \text{ and }\quad
\Ha\left(P_{\infty},Q\right)=(P_{\infty}+Q^{-1})^{-1}<P_{\infty}^{-1}.
$$
This yields the uniform estimate
\begin{eqnarray*}
\left(P_{\infty}+\widehat{\Phi}_n(Q)^{-1}\right)^{-1}&\leq& (E^n)^{\prime}\,P_{\infty}^{-1}\,
E^n+H.
\end{eqnarray*} 
Applying (\ref{wood}) to $M=H$, $N=P_{\infty}^{-1}$
and $(U,V)=((E^n)^{\prime},E^n)$ we check that
\begin{eqnarray*}
 \widehat{\Phi}_n(Q)^{-1}&\geq&\widehat{P}_{\infty}^{-1}-H^{-1}(E^n)^{\prime}\left(P_{\infty}+E^nH^{-1}(E^n)^{\prime}\right)^{-1}
E^nH^{-1}.
\end{eqnarray*}
Choosing $n_{\epsilon}$ as in (\ref{def-n-epsi}) for any $n\geq n_{\epsilon}$ we obtain the estimate
$$
\widehat{\Phi}_n(Q)^{-1}\geq \widehat{P}_{\infty}^{-1}-H^{-1}(E^n)^{\prime}P_{\infty}^{-1}
E^nH^{-1}\geq \epsilon\widehat{P}_{\infty}^{-1}.
$$
We conclude that
$$
\forall Q\in\Sa^+_r\quad \forall n\geq n_{\epsilon}\qquad \widehat{\Phi}_n(Q)\leq \epsilon^{-1}\widehat{P}_{\infty}.$$
Clearly, for any $Q\in \Sa^0_r$ by the monotone properties of Riccati semigroups stated in  lemma~\ref{lem-phi-props} for any $n\geq n_{\epsilon}$ we also have
$$
\widehat{\Phi}_n(Q)\leq \widehat{\Phi}_n(Q+I)\leq \epsilon^{-1}\widehat{P}_{\infty}
$$
This ends the proof of Corollary~\ref{theo-unif-intro}.
\cqfd

\section{Floquet-type formulae}\label{theo-2-proof}
This section is concerned with the proof of Theorem~\ref{theo-2}. 

\begin{lem}
For any $n\geq r$,  $\La_n$ maps $\Sa^0_r$ into $\Ga \La_r$. In addition, for any $n\geq r$ and $P\in \Sa^0_r$ we have the uniform estimate
\begin{equation}\label{unif-estimate}
\Vert \La_n(P)^{-1}\Vert\leq \Vert \widehat{P}_{\infty}\Vert~\Vert G_r^{-1}\Vert<\infty.
\end{equation}
\end{lem}
\proof
By Proposition~\ref{lu-estimate-Delta-prop}, the Gramian $G_n$ is invertible for any $n\geq r$. Thus, for any $n\geq r$, we have the formula
$$
\La_n(P)=\left(P+\left(G_n^{-1}-H^{-1}\right)+\left(H^{-1}-P_{\infty}\right)\right)G_n.
$$
On the other hand, we have
$$
H=\left(P_{\infty}+\widehat{P}_{\infty}^{-1}\right)^{-1}\Longrightarrow
\La_n(P)=\left(P+\left(G_n^{-1}-H^{-1}\right)+\widehat{P}_{\infty}^{-1}\right)G_n.
$$
For any $n\geq r$ and any $P\in\Sa^0_r$, (\ref{lu-estimate-Delta}) implies that
$$
0<G_r\leq G_n\leq H
\Longrightarrow
\left(P+\left(G_n^{-1}-H^{-1}\right)+\widehat{P}_{\infty}^{-1}\right)^{-1}\leq \widehat{P}_{\infty}
\quad\mbox{\rm
and}\quad
G_n^{-1}\leq G_r^{-1},
$$
which in turn implies that
$$
\Vert\La_n(P)^{-1}\Vert\leq \Vert \widehat{P}_{\infty}\Vert~\Vert G_r^{-1}\Vert<\infty,
$$
as required.\cqfd
\begin{lem}
For any $n\geq 0$ and any $P,Q\in \Sa^0_r$ we have
\begin{equation}\label{form-n}
\Ea_n(Q)=\Ea_n(P)~\left(I+(P-Q)\Ga_n(Q)\right)\quad \mbox{and}\quad \Ea_n(Q)~\La_n(Q)=E^n.
\end{equation}
\end{lem}
\proof
Again, we use induction with respect to the parameter $n$.
Since $\Ga_0(Q)=0$ and $\Ea_0(Q)=I$, the result is immediate for $n=0$. Assume that (\ref{form-n}) holds for some $n$. Then, replacing the pair $(Q,P)$ in (\ref{form-n}) by the pair $(\Phi(Q),\Phi(P))$, we check that
\begin{eqnarray*}
\Ea_{n+1}(Q)&=&\Ea_{n}(\Phi(Q))\,\Ea(Q)\\
&=&
\Ea_n(\Phi(P))~\left(I+(\Phi(P)-\Phi(Q))\Ga_n(\Phi(Q))\right)\Ea(Q).
\end{eqnarray*}
Using (\ref{form-1}) we have
\begin{eqnarray*}
\Ea_n(\Phi(P))\left(\Phi(P)-\Phi(Q)\right)\Ga_n(\Phi(Q))\Ea(Q)
&=&\Ea_{n+1}(P)(P-Q)~\Ea(Q)^{\prime}\Ga_n(\Phi(Q))\Ea(Q)\\
&=&\Ea_{n+1}(P)(P-Q)\left(\Ga_{n+1}(Q)-\Fa(Q)\right),
\end{eqnarray*}
which implies that
\begin{eqnarray*}
\Ea_{n+1}(Q)&=&\Ea_n(\Phi(P))\left(\Ea(Q)-\Ea(P)(P-Q)\Fa(Q)\right)
+\Ea_{n+1}(P)(P-Q)\Ga_{n+1}(Q)
\end{eqnarray*}
Again, using (\ref{form-1}), we conclude that
\begin{eqnarray*}
\Ea_{n+1}(Q)&=&\Ea_{n}(\Phi(P))\Ea(P)~
+\Ea_{n+1}(P)(P-Q)\Ga_{n+1}(Q)\\
&=&\Ea_{n+1}(P)~\left(I+(P-Q)\Ga_{n+1}(Q)\right),
\end{eqnarray*}
which concludes the inductive step and thus the proof.
\cqfd

\end{document}